\documentclass[11pt,leqno]{article}
\usepackage{graphicx, amsfonts, amsthm, amsxtra, amssymb, verbatim,latexsym}
\usepackage[mathscr]{euscript}
\usepackage{hyperref}
\usepackage[all,cmtip]{xy}

\begin{document}

\def\Z{\mathbb{Z}}                   
\def\Q{\mathbb{Q}}                   
\def\C{\mathbb{C}}                   
\def\N{\mathbb{N}}                   
\def\uhp{{\mathbb H}}                
\def\A{\mathbb{A}}                   
\def\dR{{\rm dR}}                    
\def\F{{\cal F}}                     
\def\Sp{{\rm Sp}}                    
\def\Gm{\mathbb{G}_m}                 
\def\Ga{\mathbb{G}_a}                 
\def\Tr{{\rm Tr}}                      
\def\tr{{{\mathsf t}{\mathsf r}}}                 
\def\spec{{\rm Spec}}            
\def\ker{{\rm ker}}              
\def\GL{{\rm GL}}                
\def\k{{\sf k}}                     
\def\ring{{ R}}                   
\def\X{{\sf X}}                      
\def\T{{\sf T}}                      
\def\Ts{{\sf S}}
\def\cmv{{\sf M}}                    
\def\BG{{\sf G}}                       
\def\podu{{\sf pd}}                   
\def\ped{{\sf U}}                    
\def\per{{\sf  P}}                   
\def\gm{{\sf  A}}                    
\def\gma{{\sf  B}}                   
\def\ben{{\sf b}}                    

\def\Rav{{\mathfrak M }}                     
\def\Ram{{\mathfrak C}}                     
\def\Rap{{\mathfrak G}}                     

\def\nov{{\sf  n}}                    
\def\mov{{\sf  m}}                    
\def\Yuk{{\sf Y}}                     
\def\Ra{{\sf R}}                      
\def\hn{{ h}}                      
\def\cpe{{\sf C}}                     
\def\g{{\sf g}}                       
\def\t{{\sf t}}                       
\def\pedo{{\sf  \Pi}}                  

\def\Der{{\rm Der}}                   
\def\MMF{{\sf MF}}                    
\def\codim{{\rm codim}}                
\def\dim{{\rm    dim}}                
\def\Lie{{\rm Lie}}                   
\def\gg{{\mathfrak g}}                

\def\u{{\sf u}}                       

\def\imh{{  \Psi}}                 
\def\imc{{  \Phi }}                  
\def\stab{{\rm Stab }}               
\def\Vec{{\rm Vec}}                 
\def\prim{{\rm prim}}                  

\def\Fg{{\sf F}}     
\def\hol{{\rm hol}}  
\def\non{{\rm non}}  
\def\alg{{\rm alg}}  

\def\bcov{{\rm \O_\T}}       

\def\leaves{{\cal L}}        

\def\GM{{\rm GM}}

\def\perr{{\sf q}}        
\def\perdo{{\cal K}}   
\def\sfl{{\mathrm F}} 
\def\sp{{\mathbb S}}  

\newcommand\diff[1]{\frac{d #1}{dz}} 
\def\End{{\rm End}}              

\def\sing{{\rm Sing}}            
\def\cha{{\rm char}}             
\def\Gal{{\rm Gal}}              
\def\jacob{{\rm jacob}}          
\def\tjurina{{\rm tjurina}}      
\newcommand\Pn[1]{\mathbb{P}^{#1}}   
\def\Ff{\mathbb{F}}                  

\def\O{{\cal O}}                     
\def\as{\mathbb{U}}                  
\def\ring{{ R}}                         
\def\R{\mathbb{R}}                   

\newcommand\ep[1]{e^{\frac{2\pi i}{#1}}}
\newcommand\HH[2]{H^{#2}(#1)}        
\def\Mat{{\rm Mat}}              
\newcommand{\mat}[4]{
     \begin{pmatrix}
            #1 & #2 \\
            #3 & #4
       \end{pmatrix}
    }                                
\newcommand{\matt}[2]{
     \begin{pmatrix}                 
            #1   \\
            #2
       \end{pmatrix}
    }
\def\cl{{\rm cl}}                

\def\hc{{\mathsf H}}                 
\def\Hb{{\cal H}}                    
\def\pese{{\sf P}}                  

\def\PP{\tilde{\cal P}}              
\def\K{{\mathbb K}}                  

\def\M{{\cal M}}
\def\RR{{\cal R}}
\newcommand\Hi[1]{\mathbb{P}^{#1}_\infty}
\def\pt{\mathbb{C}[t]}               
\def\W{{\cal W}}                     
\def\gr{{\rm Gr}}                
\def\Im{{\rm Im}}                
\def\Re{{\rm Re}}                
\def\depth{{\rm depth}}
\newcommand\SL[2]{{\rm SL}(#1, #2)}    
\newcommand\PSL[2]{{\rm PSL}(#1, #2)}  
\def\Resi{{\rm Resi}}              

\def\L{{\cal L}}                     
\def\Aut{{\rm Aut}}              
\def\any{R}                          
\newcommand\ovl[1]{\overline{#1}}    

\newcommand\mf[2]{{M}^{#1}_{#2}}     
\newcommand\mfn[2]{{\tilde M}^{#1}_{#2}}     

\newcommand\bn[2]{\binom{#1}{#2}}    
\def\ja{{\rm j}}                 
\def\Sc{\mathsf{S}}                  
\newcommand\es[1]{g_{#1}}            
\newcommand\V{{\mathsf V}}           
\newcommand\WW{{\mathsf W}}          
\newcommand\Ss{{\cal O}}             
\def\rank{{\rm rank}}                
\def\Dif{{\cal D}}                   
\def\gcd{{\rm gcd}}                  
\def\zedi{{\rm ZD}}                  
\def\BM{{\mathsf H}}                 
\def\plf{{\sf pl}}                             
\def\sgn{{\rm sgn}}                      
\def\diag{{\rm diag}}                   
\def\hodge{{\rm Hodge}}
\def\HF{{\sf F}}                                
\def\WF{{\sf W}}                               
\def\HV{{\sf HV}}                                
\def\pol{{\rm pole}}                               
\def\bafi{{\sf r}}
\def\id{{\rm id}}                               
\def\gms{{\sf M}}                           
\def\Iso{{\rm Iso}}                           

\def\hl{{\rm L}}    
\def\imF{{\rm F}}
\def\imG{{\rm G}}
\def\cy{{Calabi-Yau }}
\def\DHR{{\rm DHR }}            
\def\H{{\sf H}}            
\def\P {\mathbb{P}}                  
\def\E{{\cal E}}
\def\L{{\sf L}}             
\def\CX{{\cal X}}
\def\dt{{\sf d}}             

\newtheorem{theo}{Theorem}[section]
\newtheorem{exam}{Example}[section]
\newtheorem{coro}{Corollary}[section]
\newtheorem{defi}{Definition}[section]
\newtheorem{prob}{Problem}[section]
\newtheorem{lemm}{Lemma}[section]
\newtheorem{prop}{Proposition}[section]
\newtheorem{rem}{Remark}[section]
\newtheorem{obs}{Observation}[section]
\newtheorem{conj}{Conjecture}
\newtheorem{nota}{Notation}[section]
\newtheorem{ass}{Assumption}[section]
\newtheorem{calc}{}
\numberwithin{equation}{section}

\begin{center}
 {\LARGE\bf   Product formulas for weight two newforms}
\\

\vspace{.25in} {\large {\sc H. Movasati \footnote{Instituto de
Matem\'atica Pura e Aplicada (IMPA), Rio de Janeiro, Brazil. email:
hossein@impa.br} and Y. Nikdelan \footnote{Universidade do Estado do
Rio de Janeiro (UERJ), Instituto de Matem\'atica
e Estat\'{i}stica (IME), Departamento de An\'alise Matem\'atica, Rio de Janeiro, Brazil. e-mail: younes.nikdelan@ime.uerj.br}}} \\
\end{center}
\vspace{.25in}
\begin{abstract}
For
a weight two newform $f$ attached to an elliptic curve $E$ defined over
rational numbers we write $f=q\prod_{n=1}^\infty (1-q^n)^{g_n}, \ g_n\in\Z$ and
we observe that for some special elliptic curves $g_n$ is an increasing sequence of positive integers.
\end{abstract}
\section{Introduction}
Mathematical experiments are cheap provided that one knows which kind of experiment to do.
In this note we report such an experiment which arose from the first author's work on
quasi-modular forms, see \cite{ho14, ho14-II},  the following simple equality
\begin{equation}
\frac{q\frac{\partial}{\partial q}\eta}{\eta}=E_2
\end{equation}
where
\begin{equation}
\eta(q):= q^ \frac{1}{24}\prod_{n=1}^\infty (1-q^n),\ \ E_2(q):=\frac{1}{24}-\sum_{n=1}^\infty \frac{nq^n}{1-q^n}
\end{equation}
and K. Ono and Y. Martin's classification of weight $2$ cusp forms which are expressible using the $\eta$ function, see \cite{martin-ono}.
The most celebrated application of modular forms is the so called (arithmetic) modularity of elliptic curves.
This was originally named as  Shimura-Taniyama conjecture, and it states that for an elliptic curve $E$ defined over $\Z$, there
is a weight two newform  $f=\sum_{n=1}^\infty f_nq^n$ such that for a good prime $p$ the number of $\mathbb F_p$-rational points of $E$
is $p-f_p$.
In the case of semi-stable elliptic curves this was proved by A. Wiles in \cite{wil95}
which enabled him to complete the proof of Fermat's last theorem. For arbitrary elliptic curves
it was proved by C. Breuil, B. Conrad, F. Diamond,  and R. Taylor in  \cite{bcdt}. Tables of the pair $(E,f)$ was produced much before
the proof of arithmetic modularity theorem, see for instance the Cremona's book \cite{Cremona1992} and the webpage \cite{lmfdb}.
In this note we use the latter source and formulate a conjecture. The present text was written in 2016 and it was distributed among
few experts in the area. Since then there have been some developments and we knew of some other related works, 
 see \cite{AliMani-2016}.  However, its main conjecture is still open, and so we decided to publish it as it is.

\section{Product formulas}
Any weight two newform $f:=\sum_{i=1}^\infty f_nq^n$ can be
represented by a product formula in the form
\begin{equation}\label{eq f conj 1}
f(q)=q\prod_{n=1}^\infty (1-q^n)^{g_n},
\end{equation}
where $g_n$'s
are integer constants. This can be verified by a direct expansion of $(1-q^n)^{g_n}$ (thanks to B. Conrad for this observation). 
In \cite{martin-ono}, Martin and Ono give a
list of all weight $2$ newforms $f(q)$ that are products and quotients of the Dedekind
eta-function, i.e,
\[
f(q)=\prod_{i=1}^s\eta^{r_i}(q^{t_i})\, , \ \ s,t_1,\ldots,t_s\in \N
\ and \ r_1, r_2, \ldots, r_s \in \Z\, .
\]
Therefore, in this case $g_n$  is a repeating sequence of numbers $r_1,r_2,\ldots,r_s$. 

Our method of
computing $g_n$ (see below) is by using the logarithmic derivative of $f$ which in some sense misleading because one does not see
the integrality of $g_n$'s. Let us suppose that we can write
$f$ as in \eqref{eq f conj 1} . Then we take the logarithmic
derivative of $f$ and find
\begin{equation} \label{eq log der1}
E_f:= \frac{q\frac{\partial f}{\partial q}}{f}=q\frac{\partial
}{\partial q}\ln f=1-\sum_{n=1}^\infty g_n \cdot n \cdot
\frac{q^n}{1-q^n} \ .
\end{equation}
\href{http://w3.impa.br/~hossein/WikiHossein/files/Singular%20Codes/2018.02-Singularcode-MovasatiNikdelan-ProductFormula.txt}
{We have written a code in {\tt Singular}}, see \cite{GPS01},  which computes $g_n$'s.
This can be found in the library {\tt foliation.lib} in the first author's webpage.
The function $E_f$ is a quasi-modular form of weight two and differential order $1$,
see \cite{kaza95} and \cite{ho14}. For examples of $E,f, E_f$ see Table \ref{table1}.
These examples have a nice property that we explain it in the next section.

\section{Building blocks of weight two newforms}
\begin{conj}
\label{conj 2}
 There is an enumerable set of weight two newforms of the form  $\eta_i^{\check r_i}(q^{\check t_i})$, for some
 $\check r_i,\check t_i\in\N$, such that
 \begin{equation}\label{eq etai}
\eta_i(q)=q^{\frac{1}{\check r_i\check t_i}}\prod_{n=1}^\infty (1-q^n)^{a_{i,n}},\ \ \ a_{i,n}\in \N
\end{equation}
where $a_{i,n}$ for fixed $i$ is an increasing sequence of positive integers  with
greatest common divisor equal to one.
Further, any other weight two newform $f$
can be written as
\begin{equation}\label{eq f eta}
f(q)=\prod_{i=1}^s\eta_i^{r_i}(q^{t_i}),\ \ s\in \N,
\end{equation}
for some $r_i\in\Z$ and $t_i\in\N$.
\end{conj}
The classical Dedekind eta function is just one of the $\eta_i$ functions in \eqref{eq etai}. For this and other examples of $\eta_i$
see Table \ref{table1}.  If we apply the logarithmic derivative to the $f$ given in \eqref{eq
f eta}, then we find
\begin{equation}\label{eq logder}
\frac{q\frac{\partial f}{\partial q}}{f}=1- \sum\limits_{m =
1}^\infty  {\left( {\sum\limits_{t_i |m} {r_i \cdot a_{i,\frac{m}
{{t_i }}} } } \right)} \,m \cdot \frac{{q^m }} {{1 - q^m }}\ .
\end{equation}
Note that $f=q+O(q^2)$ and so the constant term of the equality
\eqref{eq logder} implies that:
\begin{equation}\label{19may2016-1}
\sum_{i=1}^s\frac{  r_i t_i}{\check r_i\check t_i}=1\ .
\end{equation}
Moreover, ${\rm weight}(f)=2$ and ${\rm weight}(\eta_i) \check r_i=2$ and so
\begin{equation} \label{19may2016-2}
\sum_{i=1}^s\frac{r_i}{\check r_i}=1\ .
\end{equation}
In Table \ref{table1} by
quintuple $[a_1,a_2,a_3,a_4,a_5]$ we mean the elliptic curve given
by
$$y^2+a_1xy+a_3y=x^3+a_2x^2+a_4x+a_6.$$
\begin{table}[!]
\begin{tabular}{|c|c|c|c|l|}
  \hline
  $N$ & $[a_1,a_2,a_3,a_4,a_5]$ & $\check r$ & $\check t$  & $g_n,\ \ 1\leq n \leq 12$ \\
  \hline
  36 & $[0, 0, 0, 0, 1]$ & 4 & 6  & {\tiny $1,1,1,1,1,1,1,1,1,1,1,1,\ldots$ }\\
  \hline
  37 & $[0, 0, 1, -1, 0]$ & 2 & 1   & {\tiny $1,2,3,8,16,41,97,242,598,1532,3898,10067,\ldots$ }\\
  \hline
  43 & $[0, 1, 1, 0, 0]$ & 1 & 1  & {\tiny  $2,3,4,12,22,52,114,268,608,1448,3418,8210,\ldots$} \\
  \hline
  53 & $[0, 0, 1, -1, 0]$ & 1& 1  & {\tiny  $1,3,4,7,13,31,57,123,259,559,1195,2624,\ldots$} \\
  \hline
  61 & $[1, 0, 0, -2, 1]$ & 1& 1  & {\tiny  $1,2,3,7,10,20,38,77,149,314,626,1295,\ldots$} \\
  \hline
  79 & $[1, 1, 1, -2, 0]$ & 1& 1  & {\tiny  $1,1,2,5,6,11,18,36,61,118,213,400,\ldots$} \\
  \hline
  83 & $[1, 1, 1, 1, 0]$ & 1& 1    & {\tiny  $1,1,2,4,5,11,16,31,53,97,174,330,\ldots$} \\
  \hline
  88 & $[0, 0, 0, -4, 4]$ & 1& 2    & {\tiny  $3,6,19,48,163,506,1683,5618,19123,65634,228102,797858,\ldots$} \\
  \hline
  89 & $[1, 1, 1, -1, 0]$ & 1& 1   & {\tiny  $1,1,2,3,4,10,13,25,43,79,135,246,\ldots$} \\
  \hline
  92 & $[0, 0, 0, -1, 1]$ & 1& 2  & {\tiny  $3,5,18,43,138,426,1371,4428,14683,48882,164970,560368,\ldots$} \\
  \hline
  101 & $[0, 1, 1, -1, -1]$ & 1& 1  & {\tiny  $0,2,2,2,4,7,10,18,30,52,84,152,\ldots$} \\
  \hline
  243 & $[0, 0, 1, 0, -1]$ & 1& 3  & {\tiny  $2,5,10,32,80,234,668,1988,5888,17840,54284,166950,\ldots$} \\
   & $[0, 0, 1, 0, 20]$ & & &  \\
  \hline
  256 & $[0, 0, 0, -2,0]$ & 1& 4  & {\tiny  $4,9,36,129,516,2041,8516,35780,153252,663305,2901860,12795009,\ldots$} \\
   & $[0, 0, 0, 8,0]$ &  & &  \\
  \hline
  288 & $[0, 0, 0, -12,0]$ & 2& 4  & {\tiny  $2,3,13,46,166,593,2266,8712,34147,135033,540090,2176712,\ldots$} \\
   & $[0, 0, 0,3,0]$ & &  &  \\
  \hline
  389 & $[0, 1, 1, -2, 0]$ & 1& 1 &  {\tiny  $2,3,4,11,20,51,110,259,582,1395,3262,7822,\ldots$} \\
  \hline
   675 & $[0, 0, 1, 0, -169]$ & 1& 3  & {\tiny  $2,5,10,20,56,129,362,945,2590,7093,19772,55306,\ldots$} \\
   & $[0, 0, 1, 0, 6]$ & & &  \\
  \hline
   2304 & $[0, 0, 0, -72, 0]$ & 1& 4 &  {\tiny  $4,6,16,42,132,381,1220,3851,12532,40994,135908,453455,\ldots$} \\
   & $[0, 0, 0, 18, 0]$ & &  &  \\
  \hline
\end{tabular}
\caption{Weight two newforms with a single $\eta_i$}
\label{table1}
\end{table}
A general method to verify Conjecture \ref{conj 2} is as follow. We
take two of $\eta_i$'s in Table \ref{table1}, let us say $\eta_1$
and $\eta_2$, with the corresponding conductor $N_1$ and $N_2$. Now
we search for examples of weight two newforms
$f=\eta_1^{r_1}(q^{t_1})\eta_2^{r_2}(q^{t_2})$ with explicit $r_1\in
\Z$ and $t_i\in\N$ satisfying \eqref{19may2016-1} and
\eqref{19may2016-2}. The conductor of $f$ might be $N_1N_2t_1t_2$ or
a product of powers of its primes. Similar methods as in
\cite{martin-ono} can be used in order to classify all such
newforms.

In our way to analyze Conjecture \ref{conj 2} we found some
relations with the Ramanujan theta functions. Ramanujan's
two-variable theta function $f(a,b)$ is defined by
\begin{equation}
f(a,b) = \sum_{n=-\infty}^\infty
a^{n(n+1)/2} \; b^{n(n-1)/2}
 = (-a; ab)_\infty \;(-b; ab)_\infty \;(ab;ab)_\infty,\ \ \hbox{ for   } |ab| <
 1\, ,
 \end{equation}
 where $(a,q)_n:=\prod_{i=0}^{n-1}(1-aq^i)$ is the $q$-Pochhammer symbol.
 The second identity is called  the Jacobi triple product identity.
Let us denote the newform associated with the conductor $N=256$
given in Table \ref{table1} by $f(q)=\eta_{256}(q^4)$. We get
\begin{equation}\label{eq eta256}
\eta_{256}=q^\frac{1}{4}(1-4q-3q^2-4q^3-2q^4+11q^6-4q^7+12q^9-10q^{10}+12q^{11}+7q^{12}+\ldots)\
,
\end{equation}
and hence we find
\[
\eta_{256}^2(q^2)=q-8q^3+10q^5+16q^7+37q^9+40q^{11}+50q^{13}-80q^{15}-30q^{17}-40q^{19}-128q^{21}+\ldots\,
\]
which must be an eigenform of weight $4$, because its coefficients
are multiplicative. It is in a close relationship with
\cite[$A228072$]{oeis}. We have
\[
q^{-\frac{1}{2}}\eta_{256}^2(q)=\varphi^2(q^2)\psi^2(-q^2)(\varphi^4(q^2)-8q\psi^4(-q^2))\ ,
\]
where $\varphi(q)=f(q,q)=\theta_3(q^2)$ and
$\psi(q)=f(q,q^3)=\frac{1}{2}q^{-1/4}\theta_2(q)$. We also found
that
\[
\eta_{256}^2(q)=\frac{\eta^{12}(-q^2) - 8q
\eta^{12}(q^4)}{\eta^2(-q^2)\eta^2(q^4)}.
\]
Finally, it is worth to mention that in the literature we have the  Borcherds lift, see \cite{Borcherds1995} Theorem 14.1,
in which we have 
a class of meromorphic modular forms for some character of $\SL 2\Z$, of integral weight, leading coefficient one,
whose coefficients are integers, all of whose zeros and poles are either cusps or imaginary quadratic irrationals and 
have  a product formula. The intersection of this class with the class of weight two newforms does not seem to be big, and both cases
might be generalized into a general framework in which the product formulas are explained in a uniform way.



\def\cprime{$'$} \def\cprime{$'$} \def\cprime{$'$} \def\cprime{$'$}

\end{document}